\documentclass[12pt]{extarticle}
\usepackage{latexsym}
\usepackage{amsfonts,amssymb,amsmath,mathtools}
\usepackage{bbold}

\newtheorem{thm}{Theorem}[section]

\newtheorem{lem}[thm]{Lemma}
\newtheorem{pro}[thm]{Proposition}

\newcommand{\ov }{\overline }
\newcommand{\un }{\underline }
\newcommand{\x}{\hspace{-0.025in}\times\hspace{-0.025in}}
\newcommand{\e}{\varepsilon}

\title{ The symmetric Post Correspondence Problem, and  \\
errata for the freeness problem for matrix semigroups }

\author{ J.C.\ Birget, \ A.L.\ Talambutsa }

\date{\scriptsize{
12 VII 2022}}

\begin{document}
\maketitle

\begin{abstract}
We define the {\em symmetric} Post Correspondence Problem (PCP) and prove
that it is undecidable. As an application we show that the original proof of
undecidability of the freeness problem for $3 \x 3$ integer matrix semigroups
works for the symmetric PCP, but not for the PCP in general.
\end{abstract}

\section{Introduction}

The {\em Post Correspondence Problem} (PCP) was introduced, and proved to
be undecidable, by Emil Post in 1946 \cite{Post46}.
Let $A$ be a finite alphabet of size at least 2, and let $A^*$ be the set
of all finite strings, including the empty string $\e$. The statement of
the PCP over the alphabet $A$ is as follows.

\medskip

\noindent
{\sc Input:} A non-empty finite set $\{(u_i, v_i) : 0 \le i < k\}$ of
ordered pairs of strings in $A^*$, for some $k > 0$.

\smallskip

\noindent
{\sc Question:} Does there exist a non-empty finite sequence
$(j_1,\, \ldots\, , j_n)$ of numbers in $\{0,\, \ldots\, , k-1\}$ such that
$u_{j_1}\, \ldots\, u_{j_n} = v_{j_1}\, \ldots\, v_{j_n}$?
 \ Equivalently, does the subsemigroup
$\,\langle \{(u_i, v_i) : 0 \le i < k\}\rangle\,$ of $A^* \x A^*$
intersect $\{(x,x) : x \in A^*\}$?

The PCP is called {\em bounded} if and only if some upper-bound on $k$ has been fixed
beforehand.

\medskip

\noindent Notation: $A^* \x A^*$ is the direct product of the free
monoid $A^*$ with itself. For any set $S \subseteq A^* \x A^*$,
$\,\langle S \rangle$ is the subsemigroup of $A^* \x A^*$ {\em generated}
by $S$.  Two sets $X$ and $Y$ are said to {\em intersect} if and only if
$\,X \cap Y \ne \varnothing$.
For a string $w \in A^*$, $|w|$ denotes the {\em length} of $w$; and for a
finite set $X$, $|X|$ denotes the {\em cardinality}.  We denote
{\em concatenation} of strings $u, v \in A^*$ by $uv$.
For $u, w \in A^*$ we call $u$ a {\em prefix} of $w$ if and only if there exists
$v \in A^*$ such that $w = uv$. We say that $u, v \in A^*$ are
{\em prefix-comparable} if and only if $u$ is a prefix of $v$ or $v$ is a prefix of $u$.

\medskip

\noindent {\em Non-triviality assumption:}
In order to avoid trivial solutions of the PCP (consisting of a single input
pair), we assume that the input satisfies $u_i \ne v_i$ for all $i$.

\bigskip

\noindent
The {\bf symmetric Post Correspondence Problem (symPCP)} has the same
problem statement as the PCP, but with the additional restriction that the
input relation should be {\em symmetric}; i.e., for every $(u,v) \in$
$\{(u_i, v_i) : 0 \le i < k\}\,$ we also have $\,(v,u) \in$
$\{(u_i, v_i) : 0 \le i < k\}$.

In other words, the PCP is symmetric if and only if for every $i \in \{0,\ldots,k-1\}$
there exists $i' \in \{0,\ldots,k-1\}$ such that
$\,(u_i, v_i) = (v_{i'}, u_{i'})$.

\bigskip

In Section 2 we prove that the symmetric PCP is undecidable.
In Section 3 we address some issues about the original proof of
undecidability of the freeness problem for semigroups of $3 \x 3$ integer
matrices in \cite{KlaBiSa}. In Subsection 3.2 the most significant of these
issues is resolved by reduction from the symmetric PCP instead of the
general PCP; this was our initial motivation for looking at the symmetric
PCP.

\section{Undecidability of the symmetric Post Correspondence Problem}

\begin{pro} \label{PROPsymPCPundec}
 \ For any alphabet of size at least 2, the {\em symmetric} bounded Post
Correspondence Problem is undecidable, for some bound.
\end{pro}
{\sc Proof.} Our proof is based on the proof of undecidability of the PCP by
Robert Floyd \cite{RWFloyd}, with small modifications; other references for
this proof are \cite{HdBoLogic}, \cite{DavisWeyuker}, \cite{FloydBeigel}.
Floyd's proof reduces the word problem of any semi-Thue system to a PCP.
However, since the word problem is already undecidable for certain finitely
presented semigroups (which are special semi-Thue systems, namely symmetric
Thue systems that do not use the empty string),
we immediately obtain a PCP that is almost symmetric.

More precisely, let $\langle B \mid R \rangle$ be a finite presentation of a
semigroup with undecidable word problem; here, $B$ is a finite alphabet and
$R \subseteq B^+ \x B^+$ is a finite symmetric relation. The existence of
such semigroups was proved by A.A.\ Markov and E.\ Post independently in
1947 \cite{Markov,Post47}. Note that for a semigroup presentation, only
non-empty strings are used (i.e., the set $B^+$).

For the details, we follow Floyd's proof in the formulation of
\cite{HdBoLogic, DavisWeyuker}. For the PCP we use the alphabet  $\,A \,=\,$
$\,\{\, \llcorner\,, \,\lrcorner\, , \,\circ, \, \ov{\circ}\, \}$ $\cup$
$B$  $\cup$  $\ov{B}$,  where $\ov{B} = \{\ov{b}: b \in B\}$. This alphabet
has size $|A| = 4 + 2\,|B|$, but we will later encode $A$ over $\{0,1\}$.
We will use the overline as an isomorphism from $(B \cup \{\circ\})^*$ onto
$(\ov{B} \cup \{\ov{\circ}\})^*$, with
$\,\ov{\ell_1 \ell_2 \, \ldots \, \ell_n}$  $=$
$\ov{\ell}_1 \,\ov{\ell}_2 \ \ldots \ \ov{\ell}_n \,$ for all
$\ell_1, \ell_2, \ldots, \ell_n \in B \cup \{\circ\}$.

\smallskip

An instance $\,x \stackrel{?}{=}_{_{\langle B \mid R \rangle}} y\,$ of the word
problem of the semigroup presentation $\langle B \mid R \rangle$, with
$x, y \in B^+$, is reduced to the PCP $P_{x,y}$ with the following input:

\medskip

\hspace{.2in} {\sc Input}$(P_{x,y}) \ = \ $
$\{(b, \ov{b}) : b \in B \cup \{\circ\}\}$
 \ \ $\cup$ \ \ $\{(\ov{b}, b) : b \in B \cup \{\circ\} \}$

\smallskip

\hspace{1.4in}  $\cup$ \ \ $\{(u, \ov{v}) : (u,v) \in R\}$
 \ \ $\cup$ \ \ $\{(\ov{u}, v) : (u,v) \in R\}$

\smallskip

\hspace{1.4in}  $\cup$
 \ \ $\{(\llcorner x \circ, \,\llcorner),$
$ \, (\lrcorner, \,\ov{\circ} y \lrcorner)\}\,$.

\medskip

\noindent Since $R$ can be assumed to be symmetric, we see that except for
the two pairs
$(\llcorner x\circ, \,\llcorner)$ and $(\lrcorner, \,\ov{\circ} y \lrcorner)$,
this PCP is already symmetric.
Floyd proves that this is indeed a many-one reduction, i.e.,
$x =_{_{\langle B \mid R \rangle}} y\,$ is true if and only if the PCP $P_{x,y}$ has a
solution. More precisely, there is a derivation
$\, x = x_1 \to x_2 \to \ \ldots \ \to x_{n-1} \to x_n = y\,$ in
$\langle B \mid R \rangle$ \ if and only if \ the PCP $P_{x,y}$ has a solution

\medskip

$(\llcorner x_1 \circ, \,\llcorner)$
$(\ov{x}_2, x_1) \ (\ov{\circ}, \circ)$
$(x_3, \ov{x}_2) \ (\circ, \ov{\circ})$
 \ \ $\ldots$
 \ \ $(\ov{\circ}, \circ) \ (x_n, \ov{x}_{n-1})$
$(\lrcorner, \ov{\circ} x_n \lrcorner)$

\smallskip

$=$
 \ $(\llcorner x_1\circ \ \ \ov{x}_2 \ \ \ov{\circ} \ \  x_3 \ \ \circ $
$ \ \ \ldots \ \ $   $\ov{\circ} \ \ x_n \ \ \lrcorner \ ${\bf ,}
 \ \ \ $\llcorner \ \  x_1 \ \ \circ \ \ \ov{x}_2 \ \ \ov{\circ}$
$ \ \ \ldots \ $  $\circ \ \ \ov{x}_{n-1} \ \ \ov{\circ} x_n \lrcorner)\,$.

\smallskip

\noindent Here we may assume that $n$ is odd, because the pairs
$(b,\ov{b})$ and $(\ov{b}, b)$ (for any $b \in B$) enable us to lengthen
the solution of the PCP by one block in
$\,(B^* \circ) \x (\ov{B}^* \ov{\circ})\,$ or
$\,(\ov{B}^* \ov{\circ}) \x (B^* \circ)$.

\smallskip

Note that the sets $B$ and $R$ can be kept fixed, since the given semigroup
$\langle B \mid R \rangle$ has an undecidable word problem; only $x$ and $y$
are variable in the input. Hence the $P_{x,y}$ has a bounded number of input
pairs.

\medskip

\noindent Finally, we obtain a {\sl symmetric} PCP by taking the PCP
$\,sP_{x,y}\,$ with input

\medskip

 \ \ \   \ \ \ {\sc Input}$(sP_{x,y})$  $ \ = \ $  {\sc Input}$(P_{x,y})$
 \ $\cup$ \   $\{\, (\llcorner\, , \,\llcorner x \circ),$
$ \ (\ov{\circ} y \lrcorner\,, \,\lrcorner)\,\}$.

\bigskip

\noindent {\sf Claim:} \ The symPCP $sP_{x,y}$ has a solution if and only if
the PCP $P_{x,y}$ has a solution.

\medskip

\noindent Proof of the Claim: Obviously, a solution for $P_{x,y}$ is also a
solution for $sP_{x,y}$.

Conversely, suppose $sP_{x,y}$ has a solution.
This solution starts either with the pair
$\,(\llcorner x\circ, \,\llcorner)\,$ or the pair
$\,(\llcorner, \,\llcorner x\circ)$, since those are the only
pairs in $sP_{x,y}$ in which the two coordinates have a common prefix;
in all other pairs, one coordinate starts with an overlined letter and
the other coordinate starts with a non-overlined letter.

Case 1: The start pair is $(\llcorner x\circ, \,\llcorner)$.
Now by the same reasoning as in \cite[pp.\ 131-132]{DavisWeyuker} and
\cite{HdBoLogic}, we can construct a derivation
$x =_{_{\langle B \mid R \rangle}} y$. In this construction, the 1st coordinate
is always longer than the 2nd coordinate, until the derivation of $y$ is
complete; then $(\lrcorner, \,\ov{\circ} y \lrcorner)$ is the right-most
pair of the solution of the PCP.

Case 2: The start pair is $(\llcorner, \,\llcorner x\circ)$.
Then by just switching the roles of the 1st and 2nd coordinates, we can
carry out the same reasoning as in Case 1; now
$(\ov{\circ} y \lrcorner\,, \,\lrcorner)$ is the right-most pair of the
solution of the PCP.
Again a derivation $x =_{_{\langle B \mid R \rangle}} y$ is constructed.

[This proves the Claim.]

\medskip

We still have to show that symPCP is undecidable for an alphabet of size 2,
e.g., for $\{0,1\}$. The symPCP $sP_{x,y}$ uses the alphabet
$\,A \,=\,$
$\{\llcorner\,, \,\lrcorner\, , \,\circ, \, \ov{\circ} \}$ $\cup$  $B$
$\cup$  $\ov{B}\,$,  of size $\,|A| = 4 + 2\,|B|$.
Let us choose any injective function
$ \ {\sf code}: A \,\to\, \{0,1\}^{\ell}$, where
$\ell = \lceil \log_2 |A| \rceil$.
Let $ \ {\sf code}(sP_{x,y})$ $=$
$\{({\sf code}(q), {\sf code}(r)) : (q,r) \in$ {\sc Input}$(sP_{x,y})\}$,
and $ \ {\sf code}(P_{x,y})$ $=$
$\{({\sf code}(q), {\sf code}(r)) : (q,r) \in$ {\sc Input}$(P_{x,y})\}$.
Obviously, ${\sf code}(sP_{x,y})$ is a symmetric PCP.

Then the symPCP $\,{\sf code}(sP_{x,y})\,$ has a solution if and only if the original
symPCP $sP_{x,y}$ has a solution.
 \hfill $\Box$

\section{Clarifications and errata for the freeness problem }

Let ${\mathbb N}^{3 \times 3}$ denote the monoid of 3-by-3 matrices
over the natural numbers. For a subset $S$ of ${\mathbb N}^{3 \times 3}$,
the subsemigroup generated by $S$ in ${\mathbb N}^{3 \times 3}$ is denoted
by $\langle S \rangle$.
Article \cite{KlaBiSa} considers  the following problem, called the
{\em freeness problem} of subsemigroups of ${\mathbb N}^{3 \times 3}$.

\smallskip

{\sc Input:} A finite set $S \subseteq {\mathbb N}^{3 \times 3}$.

\smallskip

{\sc Question:} Is $\langle S \rangle$ free over $S$?
 \ \ (Note that this is not equivalent to just asking whether
$\langle S \rangle$ is free, i.e., is isomorphic to any free semigroup.)

\smallskip

\noindent The freeness problem is shown to be undecidable in
\cite{KlaBiSa}, but the proof is incomplete (see Subsection 3.2 below, where
the gap is filled).
Many stronger forms of this result were proven later; e.g., the problem is
undecidable for upper triangular matrices in
${\mathbb N}^{3 \times 3}\,$ \cite{CaHaKa} (whose proof is not  based on
\cite{KlaBiSa}).

\bigskip

Subsection 3.1 corrects notational errors in \cite{KlaBiSa}, arising from
mix-ups between reverse base 2 and ordinary base 4 notations.

In Subsection 3.2, the main claim of \cite{KlaBiSa} is proved by using the
symmetric PCP. In \cite{KlaBiSa} it is claimed that certain finitely
generated matrix semigroups have a relation if and only if the PCPs encoded by these
matrices have solutions; this claim is true if the PCPs are symmetric.
For non-symmetric PCPs there are counter-examples, found by the second author.

\subsection{Encoding a PCP by matrices}

\noindent $\bullet$ {\it Page 224, bottom paragraph of
{\rm \cite{KlaBiSa}}:} \ The word ``reverse'' should be removed.

\medskip

\noindent {\em Comment on this correction:}
The function $\varphi: \{0,1,2,3\}^* \to {\mathbb N}\,$ performs base 4
conversion, i.e., $\varphi(x_{n-1} \ldots x_1 x_0)$  $\,=\,$
$\sum_{i=0}^{n-1} x_i \, 4^i\,$; this is the ordinary base 4 representation,
not reverse base 4.
On the other hand, the function $\,\beta: \{0,1\}^* \to {\mathbb N}\,$
performs reverse base 2 conversion, i.e.,
$\beta(x_0 x_1 \ldots x_{n-1})$ $\,=\, $ $\sum_{i=0}^{n-1} x_i \, 2^i$.

The reason for the difference is that $\varphi$ is used in a lower-triangular
$2 \x 2$ matrix, whereas $\beta$ is used in an upper-triangular
$2 \x 2$ matrix.

\bigskip

\noindent $\bullet$ {\it The first paragraph on page 225 of
{\rm \cite{KlaBiSa}} should be replaced by the following:}

\smallskip

Next we want to encode an instance of the Post Correspondence Problem into
$3 \x 3$ matrices over ${\mathbb N}$.
We view the indices $\,0, 1, \ldots, k-1\,$ of the Post Correspondence
pairs as binary strings of uniform length
$\,h = \max\{\lceil \log_2 k \rceil,\, 1\}$.
More precisely, we encode every $\,i \in \{0, 1, \ldots, k-1\}\,$ as
$\,\un{i} \in \{0, 1\}^h$; the choice is arbitrary, except that
$i \mapsto \un{i}$ is injective.
Since $\{0,1\} \subseteq \{0, 1, 2, 3\}$ we can decode the string $\un{i}$
into a natural number by $\varphi$, in ordinary base 4 notation.
The binary strings $u_i$ and $v_i$ will be decoded into natural numbers by
$\beta(u_i)$, respectively $\beta(v_i)$, in reverse binary notation.

Then we represent the instance $\,\{(u_i, v_i) :$
$i = 0,1,\, \ldots\, , k-1\}\,$ of the Post Correspondence Problem by
$\,4k+1\,$ matrices as follows:

$$L = \left[
\begin{array}{ccc}
1 & 0 & 0 \\
0 & 1 & 0 \\
0 & 2 & 4
\end{array}  \right] \ , $$
and for $\,i = 0, 1,\, \ldots\, , k-1$:
$$ U_i \ = \ \left[
\begin{array}{ccc}
2^{|u_i|} & \beta(u_i)   & 0 \\
0         & 1            & 0 \\
0         & 2 + \varphi(\un{i}) \cdot 4  & 4^{h+1}
\end{array}       \right]\, ,
$$

$$ \ov{U}_i \ = \ \left[
\begin{array}{ccc}
2^{|u_i|} & \beta(u_i)   & 0 \\
0         & 1            & 0 \\
0         & 2 + \varphi(\un{i}) \cdot 4 + 3 \cdot 4^{h+1} & 4^{h+2}
\end{array}       \right]\, ,
$$

$$ V_i \ = \ \left[
\begin{array}{ccc}
2^{|v_i|} & \beta(v_i)   & 0 \\
0         & 1            & 0 \\
0         & \varphi(\un{i}) + 2 \cdot 4^h  & 4^{h+1}
\end{array} \right]\, ,
$$

$$ \ov{V}_i \ = \ \left[
\begin{array}{ccc}
2^{|v_i|} & \beta(v_i)   & 0 \\
0         & 1            & 0 \\
0         & \varphi(\un{i}) + 3 \cdot 4^h  & 4^{h+1}
\end{array} \right]\, .
$$

\bigskip

\noindent In summary, the correspondence between these matrices and pairs
of strings is given by the following table (where
$i = 0, 1, \,\ldots\,, k-1$):

\bigskip

\begin{table}
\hspace{1.5in}  \begin{tabular} {|l|l|l|l|l|}
$L$      & $U_i$             & $\ov{U}_i$          & $V_i$           &          $\ov{V}_i$ \\ \hline
$(\e,2)$ & $(u_i, 2 \un{i})$ & $(u_i, 2 \un{i} 3)$ & $(v_i, \un{i} 2)$ & $(v_i, \un{i} 3)$
\end{tabular}
\caption{\label{CorrTable}Correspondence between matrices and pairs of strings.}
\end{table}

\bigskip

\noindent As we remarked earlier, when these matrices are multiplied, the
pairs of strings that they encode are concatenated; for example,
$U_iU_j\ov{U}_r$ encodes $(u_iu_ju_r,\, 2\un{i} 2 \un{j} 2 \un{r}3)$.
Note again that all strings $\un{i}$ are over $\{0,1\}$ and have
the same length $h = \max\{\lceil \log_2 k \rceil,\, 1\}$.

\medskip

\noindent {\em Comment about this correction:}
In \cite{KlaBiSa}, $i$ appears in some places where $\varphi(\un{i})$
should have been used (although the explanations in the paper make it clear
that $\varphi(\un{i})$ was intended).
However, $\varphi(\un{i})$ is a different number than $i$. Indeed, $\un{i}$
is an arbitrary binary string representing $i$, whereas $\varphi(\un{i})$
is the natural integer represented by the binary string $\un{i}$ in base 4;
note also that this is the ordinary base 4 decoding, not reverse base 4.

\bigskip

\noindent $\bullet$ {\it Isomorphism between subsemigroups of
 ${\mathbb N}^{3 \times 3}$ and subsemigroups of
$\,\{0,1\}^* \x \{0,1,2,3\}^*$}

\medskip

\noindent The following was briefly mentioned in \cite{KlaBiSa} (p.\ 224,
bottom paragraph); here we give more details:

\begin{lem} \label{LEMembtodirect}
 \ The subsemigroup of ${\mathbb N}^{3 \times 3}$ generated by
$\,\{L\} \,\cup\, \{U_i, \ov{U}_i, V_i, \ov{V}_i : 0 \le i < k\}\,$ can be
embedded into $\,\{0,1\}^* \x \{0,1,2,3\}^* \ $ (direct product of two free
monoids).  More precisely,

\medskip

 \ \ \ $\big\langle \{L\}$ $\,\cup\,$
$\{U_i, \ov{U}_i, V_i, \ov{V}_i : 0 \le i < k\} \big\rangle \ $ in
$\,{\mathbb N}^{3 \times 3}$

\medskip

\noindent is isomorphic to

\medskip

 \ \ \ $\big\langle \{(\e,2)\}$  $\,\cup\,$
$\{(u_i, 2\un{i}), (u_i, 2\un{i}3), (v_i, \un{i}2), (v_i, \un{i}3):$
 $0 \le i < k\} \big\rangle \ $ in $\,\{0,1\}^* \x \{0,1,2,3\}^*$

\medskip

\noindent by the isomorphism given by the Table~\ref{CorrTable} above.
\end{lem}
{\sc Proof.}  As we indicated above, the matrices represent
pairs of strings, and we remarked that when these matrices are multiplied,
the pairs of strings that they encode are concatenated.

For any matrix $M$, let $\,M_{i,j}$ denote the $(i,j)$-entry (in row $i$
and column $j$).

\medskip

\noindent {\sf Claim:} Let $M$ be a product of a sequence of matrices in
$\{L\}$ $\cup$ $\{U_i, \ov{U}_i, V_i, \ov{V}_i : 0 \le i < k\}$.
Then the entries $\,M_{1,1}, M_{1,2}, M_{3,3}, M_{3,2}\,$ determine a pair
$(w, J) \in$ $\{0,1\}^* \x \{0,1,2,3\}^*$, where $(w, J)$ is the
concatenation of the corresponding pairs in
$\,\{(\e,2)\}$  $\cup$
$\{(u_i, 2\un{i}), (u_i, 2\un{i}3), (v_i, \un{i}2), (v_i, \un{i}3):$
$0 \le i < k\}$.
Conversely, $(w, J)$ determines $M$ by
$$(w, J) \ \ \ \longmapsto \ \ \   M \,=\, \left[
\begin{array}{ccc}
2^{|w|}   & \beta(w)        & 0 \\
0         & 1               & 0 \\
0         & \varphi(J) & 4^{|J|}
\end{array}       \right]\,.
$$

\medskip

\noindent Note that $w$ $\in \{0,1\}^*$ could have trailing 0s, so in
reverse base 2 representation, $\beta(w)$ $\in {\mathbb N}$ alone does
not determine $w$. But $\beta(w)$ and $|w|$ together determine $w$.
Similarly, $J$ $\in \{0,1,2,3\}^*$ could have leading 0s, so in base 4
representation, $\varphi(J)$ $\in {\mathbb N}$ alone does not determine
$J$. But $\varphi(J)$ and $|J|$ together determine $J$.

\smallskip

\noindent Proof of the Claim: We use induction on the number of matrices
multiplied.
For one matrix, the lemma holds by the Table~\ref{CorrTable} above.
In general, suppose $M$ determines $(w,J)$, as in the lemma,
and consider the matrix $MX$, where $\,X \in \{L\}$ $\cup$
$\{U_i, \ov{U}_i, V_i, \ov{V}_i : 0 \le i < k\}$.

\medskip

\noindent If $X = L$,

\smallskip

$ ML \ = \ \left[
\begin{array}{ccc}
2^{|w|}   & \beta(w)        & 0 \\
0         & 1               & 0 \\
0         & \varphi(J) + 2 \cdot 4^{|J|} & 4^{|J|+1}
\end{array}       \right]
 \ = \ \left[
\begin{array}{ccc}
2^{|w|}   & \beta(w)        & 0 \\
0         & 1               & 0 \\
0         & \varphi(J2) & 4^{|J2|}
\end{array}       \right] \,,
$

\medskip

\noindent which determines $(w, J2)$. Recall that the 1st coordinate of
$(w,J)$ uses reverse base 2, and the 2nd coordinate uses the usual base 4
representation.

\medskip

\noindent If $X = U_i$,

\smallskip

$MU_i \ = \ \left[
\begin{array}{ccc}
2^{|w| + |u_i|}   & 2^{|w|}\,\beta(u_i) + \beta(w) & 0 \\
0                 & 1                                     & 0 \\
0                 & \varphi(J) + (2 + \varphi(\un{i}) \cdot 4) \cdot 4^{|J|}  & 4^{|J|+h+1}
\end{array}   \right]
 \ = \ \left[
\begin{array}{ccc}
2^{|w u_i|}   & \beta(w u_i)  & 0 \\
0                 & 1                   & 0 \\
0                 & \varphi(J 2 \un{i}) & 4^{|J 2 \un{i}|}
\end{array}       \right]\,,
$

\medskip

\noindent which determines $(wu_i, J 2 \un{i})$.

\medskip

\noindent If $X = \ov{U}_i$,

\smallskip

$M \ov{U}_i \ = \ \left[
\begin{array}{ccc}
2^{|w| + |u_i|} & 2^{|w|}\,\beta(u_i) + \beta(w) & 0 \\
0               & 1                             & 0 \\
0 & \varphi(J) + (2+\varphi(\un{i})\cdot4+3\cdot 4^{h+1})\cdot4^{|J|} &
4^{|J|+h+2}
\end{array}       \right]$

\smallskip

\hspace{.38in} $= \ \left[
\begin{array}{ccc}
2^{|w u_i|}   & \beta(w u_i)  & 0 \\
0                 & 1                   & 0 \\
0                 & \varphi(J 2 \un{i} 3) & 4^{|J 2 \un{i} 3|}
\end{array}       \right],
$

\medskip

\noindent which determines $(wu_i, J 2 \un{i} 3)$.

\smallskip

Similarly, $M V_i$ and $M \ov{V}_i$ determine $(w, J \un{i} 2)$,
respectively $(w, J \un{i} 3)$.

[This proves the Claim.]

\medskip

\noindent The function from matrices to pairs of strings is injective, since
$(w,J)$ determines $M$ by the formula given in the claim.
\hfill $\Box$

\bigskip

\noindent As a consequence of Lemma \ref{LEMembtodirect} we have:

\medskip

$\langle \{L\}$ $\,\cup\,$
$\{U_i, \ov{U}_i, V_i, \ov{V}_i : 0 \le i < k\} \rangle\,$ is free over the
given generators

\medskip

if and only if \ \ $\langle \{(\e,2)\}$  $\,\cup\,$
$\{(u_i, 2\un{i}), (u_i, 2\un{i}3), (v_i, \un{i}2), (v_i, \un{i}3):$
 $0 \le i < k\} \rangle\,$ is free over the given generators.

\subsection{An application of the symmetric PCP}

The proof in \cite{KlaBiSa} does not work for all PCPs. A counter-example was
found by A.\ Talambutsa, who also observed that restricting the proof to the
symmetric PCP would correct the mistake.

\bigskip

\noindent {\bf Example:} Consider the PCP with input $\{(00,0)\}$, which
obviously has no solution. (There are similar counter-examples, e.g.,
$\{(00,0),\, (u,0)\}$, for any $u \in \{0,1\}^+$ with $|u| \ge 2$.)
The corresponding subsemigroup of $\{0,1\}^* \x \{0,1,2,3\}^*$, constructed
in \cite{KlaBiSa}, is generated by $ \ \Gamma$  $=$
$\{(\e,2),$  $(00, 20),$  $(00, 203),$  $(0,02),$  $(0,03)\}$;
here the single input pair $(00,0)$ is coded by the binary string
$\un{0} = 0\,$ (as in Subsection 3.1).
The subsemigroup $\langle \Gamma \rangle$ satisfies the following relations
(among others):

\smallskip

 \ \ \    $(00, 20)\,(\e,2)\,(\e,2)\,(0,03) $
$ \ = \ $ $(\e,2)\,(0,02)\,(00, 203) \,$;

\smallskip

 \ \ \    $(00, 20)\,(\e,2)\,(\e,2)\,(0,02)$
$ \ = \ $  $(\e,2)\,(0,02)\,(00, 20)\,(\e,2) \,$;

\smallskip

 \ \ \    $(00, 20)\,(00, 203)\,(\e,2)\,(0, 02)\,(0,03) $
$ \ = \ $ $(\e,2)\,(0,02)\,(0,03)\,(00, 20)\,(00, 203) \,$.

\smallskip

\noindent So, the semigroup has relations that do not correspond to
solutions of the PCP. The 1st coordinate of these relations does, however,
correspond to solutions of the {\sl symmetric} PCP $\,\{(00,0), (0,00)\}$.

\bigskip

\noindent Next we show in detail that the proof in \cite{KlaBiSa}, which
was not very detailed, works correctly for the symmetric PCP.
From here on we consider the input of the PCP to be a {\sl sequence} of
different pairs, rather than a set; let
$\,\big( (u_i,v_i) : 0 \le i < k \big)\,$ be a sequence corresponding to
some chosen ordering of the pairs of the PCP input
$\{(u_i,v_i) : 0 \le i < k\}$.

\begin{lem} \label{LEMreduction} {\bf (reduction).}
 \ The {\em symPCP} with input $\big((u_i,v_i) : 0 \le i < k\big)$ over
$\{0,1\}$ has a solution \ if and only if \ the subsemigroup of
$\,\{0,1\}^* \x \{0,1,2,3\}^*\,$ generated by

\medskip

 \ \ \  $\Gamma \ = \  \{(\e,2)\}$  $ \ \cup \ $
$\{(u_i, 2 \un{i}),\,(u_i, 2 \un{i} 3),\,(v_i, \un{i} 2),\,(v_i, \un{i} 3)$
$\,:\, 0 \le i < k\}$

\medskip

\noindent is not free. The function
$\,\big((u_i,v_i) : 0 \le i < k\big)$ $\,\longmapsto\,$
$\Gamma\,$ is a one-one polynomial-time reduction from
the symPCP to the non-freeness problem of the subsemigroup
$\langle \Gamma \rangle$ of $\,\{0,1\}^* \x \{0,1,2,3\}^*$.
\end{lem}
{\sc Proof.} $[\Rightarrow]$ \ If the PCP has a solution
$ \ u_{j_1} \ \ldots \ u_{j_n} = v_{j_1} \ \ldots \ v_{j_n} \ $
then one verifies immediately  that the following semigroup relation holds
in $\langle \Gamma \rangle$:

\smallskip

 \ \ $(u_{j_1}, 2 \un{j_1})$   $ \ \ldots \ $
$(u_{j_{n-1}}, 2 \un{j_{n-1}}) (u_{j_n}, 2 \un{j_n} 3)$
$ \ = \ $
$(\e,2) (v_{j_1}, \un{j_1} 2)$  $ \ \ldots \ $
$(v_{j_{n-1}}, \un{j_{n-1}} 2) (v_{j_n}, \un{j_n} 3)$.

\medskip

\noindent So $\langle \Gamma \rangle$ is not free.

Also, using symmetry we have $(u_{j_n}, v_{j_n})$  $=$
$(v_{j'_n}, u_{j'_n}\,)$ for some $j'_n \in \{0,\ldots,k-1\}$. So
$\,(u_{j_n}, \un{j'_n} 2)$, $(v_{j_n}, 2 \un{j'_n})$  $\in \Gamma\,$ \big(in
addition to $\,(u_{j_n}, 2 \un{j_n})$, $(v_{j_n}, \un{j_n} 2)$
$\in \Gamma\,$\big).
Note that $(u_i,v_i)$ has a code $\un{i} \in \{0,1\}^*$ which is different
from the code $\un{i'}$ of $(v_i,u_i) = (u_{i'}, v_{i'})$ \ (by the
non-triviality assumption on {\sc Input}(PCP)).
 \ Hence we also have the relation,

\medskip

$(u_{j_1}, 2 \un{j_1})$  $ \ \ldots \ $
$(u_{j_{n-1}}, 2 \un{j_{n-1}}) (\e,2)(\e,2) (u_{j_n}, \un{j'_n} 2)$
$ \ = \ $
$(\e,2) (v_{j_1}, \un{j_1} 2)$   $ \ \ldots \ $
$(v_{j_{n-1}}, \un{j_{n-1}} 2) (v_{j_n}, 2 \un{j'_n}) (\e, 2)$.

\medskip

\noindent The corresponding matrix relations are

\smallskip

 \ \ \ $U_{j_1} \ \ldots \ U_{j_{n-1}} \ov{U}_{j_n}$
$\,=\,$   $L V_{j_1} \ \ldots \ V_{j_{n-1}} \ov{V}_{j_n}$,
 \ \ \  \ \ and

\smallskip

 \ \ \ $U_{j_1} \ \ldots \ U_{j_{n-1}} L^2\, V_{j'_n}$
$\,=\,$   $L V_{j_1} \ \ldots \ V_{j_{n-1}} U_{j'_n} L$.

\medskip

\noindent The first of these matrix relations was given in \cite{KlaBiSa},
the second is new and is based on the symmetry of the PCP input.

\bigskip

\noindent {\bf Remark (\ref{LEMreduction}R):}
The second, new, relation (based on the symPCP) shows that the letter 3,
as well as the generators
$\,\{(u_i, 2 \un{i} 3),\, (v_i, \un{i} 3)$ $: 0 \le i < k\}$,
are not needed for the construction of a subsemigroup of
$\{0,1\}^* \x \{0,1,2\}^*$ with undecidable freeness problem.
Similarly, the overlined matrices
$\,\{\ov{U}_i,\,\ov{V}_i : 0 \le i < k\}\,$
are not needed in the input of the freeness problem for matrices.

We will nevertheless continue using 3 and the redundant generators and
matrices, since we want to show that the proof in \cite{KlaBiSa} is correct
for symmetric PCPs. It is straightforward to rewrite the proof without
the redundancies (by simply leaving out the redundancies); see Lemma
\ref{LEMreduction2}. \ \ [End, Remark.]

\bigskip

\noindent $[\Leftarrow]$ \ Suppose the semigroup $\langle \Gamma \rangle$ is
not free, i.e., it has a non-trivial relation
$p_1 \,\ldots\, p_m = q_1 \,\ldots\, q_n$,
with $p_1,$ $\ldots,$ $p_m,$ $q_1,$ $\ldots,$ $q_n$  $\in$  $\Gamma$.

Since the semigroup $\{0,1\}^* \x \{0,1,2,3\}^*$ is cancellative,
we can assume that $p_1 \ne q_1$ and $p_m \ne q_n$.
We abbreviate the generator sequence $(p_1, \,\ldots\,, p_m)$ by $P$, and
the generator sequence $(q_1, \,\ldots\,, q_n)$ by $Q$. The relation is
abbreviated by $\Pi P = \Pi Q$, where $\Pi P$ is the product, in
$\{0,1\}^* \x \{0,1,2,3\}^*$, of the generators in the sequence $P$; and
similarly for $\Pi Q$.

\bigskip

\noindent {\bf Claim 1:} \\
{\sf (A)} One of the generators $p_1, q_1$ is $(\e,2)$, and the
other belongs to $\{(u_i, 2 \un{i}),\, (u_i, 2 \un{i} 3) : 0 \le i < k\}$.

\smallskip

\noindent {\sf (B)} If $q_1 = (\e,2)$, we have:

\smallskip

 \ \ \ \ $p_1 = (u_{j_1}, 2 \un{i})\,$ and
$\,q_1 q_2 = (\e,2) (v_{j_1}, \un{i} 2) = (v_{j_1}, 2 \un{i} 2)$,
 \ for some $\,(u_{j_1},v_{j_1})$  $\in$ {\sc Input}(PCP);

 \ \ \ \ or

 \ \ \ \ $p_1 = (u_{j_1}, 2 \un{i} 3)\,$ and
$\,q_1 q_2 = (\e,2) (v_{j_1}, \un{i} 3) = (v_{j_1}, 2 \un{i} 3)$,
 \ for some $\,(u_{j_1},v_{j_1})$  $\in$ {\sc Input}(PCP).

\smallskip

 \ If $p_1 = (\e,2)$ then, symmetrically, the conclusion is similar.

\medskip

\noindent Proof of Claim 1(A):  There are several cases.
By cancellativity we already ruled out $p_1 = q_1$.

\smallskip

\noindent Case (1): \ $p_1 = (\e,2)$.

Then $q_1 \in $
$\{(u_i, 2 \un{i}),\, (u_i, 2 \un{i} 3) : 0 \le i < k\}$, since the 2nd
coordinate of $q_1$ must start with 2, and we ruled out $p_1 = q_1$.

%

\smallskip

\noindent Case (2):
 \ $p_1 \in \{(u_i, 2 \un{i}),\, (u_i, 2 \un{i} 3) : 0 \le i < k\}$.

Then $q_1 \in \{(\e,2)\}$  $\cup$
$\{(u_i, 2 \un{i}),\, (u_i, 2 \un{i} 3) : 0 \le i < k\}$, as the 2nd
coordinate of $q_1$ must start with 2.

If $p_1 = (u_i, 2 \un{i})$ and $q_1 = (u_j, 2 \un{j})$ then
$\un{i}$ and $\un{j}$ are prefix-comparable, hence $\un{i} = \un{j}$
since both have length $\lceil \log_2 k \rceil$; this implies $p_1 = q_1$, which was ruled out.
Similarly, if $p_1 = (u_i, 2 \un{i} 3)$ and $q_1 = (u_j, 2 \un{j} 3)$ then
$p_1 = q_1$, which is ruled out.

If $p_1 = (u_i, 2 \un{i})$ and $q_1 = (u_j, 2 \un{j} 3)$ then
$\un{i}$ and $\un{j}$ are prefix-comparable, hence $\un{i} = \un{j}$
since both have length $\lceil \log_2 k \rceil$.
Then there will be no
possible choice for $p_2$ that could match the letter 3 in $q_1$; this
contradicts the assumption that $p_1 \ldots p_m = q_1 \ldots q_n$.

If $p_1 = (u_j, 2 \un{j} 3)$ and $q_1 = (u_i, 2 \un{i})$ or
$q_1 = (u_j, 2 \un{j} 3)$, then we obtain the same contradictions as above,
with the roles of $p_1$ and $q_1$ switched.

The only alternative left is $q_1 = (\e, 2)$.

\smallskip

\noindent Case (3):
$p_1 \in \{(v_i, \un{i} 2),\, (v_i, \un{i} 3) : 0 \le i < k\}$.

By symmetry of the input PCP, this is the same as Case (2).

\medskip

\noindent Proof of Claim 1(B):  By Claim 1(A) we can assume that
$q_1 = (\e,2)$ (the case where $p_1 = (\e,2)$ is similar). Moreover,
$p_1 = (u_{j_1}, 2 \un{j_1})$ or $p_1 = (u_{j_1}, 2 \un{j_1} 3)$, for some
$j_1$ (uniquely determined by $\Pi P$ and $\Pi Q$).
The fact that all $\un{j_i}$ have the same length $\lceil \log_2 k \rceil$
implies:
$\Pi Q = q_1q_2 \ldots$ \ for some $q_2 \in \Gamma$, and
$q_2 = (v_{j_1}, \un{j_1} 2)$, or
$q_2 = (v_{j_1}, \un{j_1} 3)$ with
$(u_{j_1},v_{j_1}) \in$ {\sc Input}(PCP).

\smallskip

[This proves Claim 1.]

\bigskip

\noindent {\bf Claim 2:}
The relation $\Pi P = \Pi Q$ has one of two forms:

\smallskip

\noindent [2-2 block]: \ \ \
$\Pi P \ = \ $
$ (u_{j_1}\,\ldots\,u_{j_r}, \ 2 \un{j_1}\,\ldots\,2 \un{j_r} 2) \ R_1$
$ \ = \ $
$\,(v_{j_1}\,\ldots\,v_{j_r}, \ 2 \un{j_1}\,\ldots\,2 \un{j_r} 2) \ R_2$
$ \ = \ \Pi Q$,

\smallskip

\noindent [2-3 block]: \ \ \
$\Pi P \ = \ $
$(u_{j_1}\,\ldots\,u_{j_r}, \ 2\un{j_1}\,\ldots\,2\un{j_r}3) \ R_3$
$ \ = \ $
$(v_{j_1}\,\ldots\,v_{j_r},\,2\un{j_1}\,\ldots\,2\un{j_r}3) \ R_4$
$ \ = \ \Pi Q$,

\smallskip

\hspace{0.67in} for some $R_1, R_2, R_3, R_4 \in \langle \Gamma \rangle$.

\smallskip

\noindent In either case we have:
 \ $(u_{j_1},v_{j_1}), \,\ldots\, ,(u_{j_r},v_{j_r})$  $\in$
{\sc Input}(PCP), and
$\,u_{j_1}\,\ldots\,u_{j_r}$ is {\sl prefix-comparable} with
$\,v_{j_1}\,\ldots\,v_{j_r}$. (Note that $\,u_{j_1}\,\ldots\,u_{j_r}\,$ and
$\,v_{j_1}\,\ldots\,v_{j_r}\,$ need not be equal.)

\smallskip

\noindent Proof of Claim 2:
By Claim 1(B) we can have

 \ \ \ $\,\Pi P = p_1\,\ldots\,$ $=$ $(u_{j_1}, 2 \un{j_1} 3)\,\ldots\,$
$=$ $\,(v_{j_1}, 2 \un{j_1} 3) \,\ldots\,$  $=$ $q_1q_2 \,\ldots\,$
$=$ $\Pi Q$,

\noindent where $(u_{j_1}, v_{j_1}) \in$ {\sc Input}(PCP), and $u_{j_1}$ is
prefix-comparable with $v_{j_1}$.
Then $\Pi P = \Pi Q$ begins with a 2-3 block, so Claim 2 holds.

\smallskip

\noindent Or we have

 \ \ \ $\,\Pi P = p_1\,\ldots\,$ $=$ $(u_{j_1}, 2 \un{j_1})\,\ldots\,$
$=$ $\,(v_{j_1}, 2 \un{j_1} 2) \,\ldots\,$  $=$ $q_1q_2 \,\ldots\,$
$=$ $\Pi Q$,

\noindent where $(u_{j_1}, v_{j_1}) \in$ {\sc Input}(PCP), and $u_{j_1}$ is
prefix-comparable with $v_{j_1}$.

\smallskip

The presence of two letters 2 in the 2nd coordinate of $q_1 q_2$ implies
that $\Pi P = p_1 p_2 \ldots \ $, for some $p_2 \in \Gamma$ of the form
$p_2 = (\e,2)$, or $p_2 = (u_{j_2}, 2 \un{j_2})$, or
$p_2 = (u_{j_2}, 2 \un{j_2} 3)$.
Hence the relation $\Pi P = \Pi Q$ takes one of the following forms:

\smallskip

\noindent Case 1: \ \ \
$\,(u_{j_1}, 2 \un{j_1} 2) \,\ldots\,$ $=$
$(v_{j_1}, 2 \un{j_1} 2) \,\ldots \ $, with $p_2 = (\e,2)$.

\smallskip

\noindent In this case the relation $\Pi P = \Pi Q$ begins with a 2-2 block,
 so Claim 2 holds.

\smallskip

\noindent Case 2: \ \ \
$\,(u_{j_1} u_{j_2}, 2 \un{j_1} 2 \un{j_2} 3) \,\ldots \ $ $=$
$(v_{j_1}, 2 \un{j_1} 2) \,\ldots \ $, \ where
$p_2 = (u_{j_2}, 2 \un{j_2} 3)$.

\smallskip

\noindent Now $\Pi Q = q_1q_2q_3\ldots \ $, for some $q_3$ of the form
$q_3 = (v_{j_2}, \un{j_2} 3) \in \Gamma$. Then the relation
$\,\Pi P= p_1 p_2 \,\ldots\, = q_1 q_2 q_3 \,\ldots\, = \Pi Q$ takes the form

\smallskip

 \ \ \ $(u_{j_1} u_{j_2},\, 2 \un{j_1} 2 \un{j_2} 3) \,\ldots \ $ $=$
   $(v_{j_1} v_{j_2},\, 2 \un{j_1} 2 \un{j_2} 3) \,\ldots \ $,

\smallskip

\noindent where $(u_{j_1}, v_{j_1}), (u_{j_2}, v_{j_2}) \in$
{\sc Input}(PCP), and $u_{j_1} u_{j_2}$ is prefix-comparable with
$v_{j_1} v_{j_2}$.
In this case the relation $\Pi P = \Pi Q$ begins with a 2-3 block, so
Claim 2 holds.

\smallskip

\noindent Case 3: \ \ \
$\,(u_{j_1} u_{j_2}, 2 \un{j_1} 2 \un{j_2}) \,\ldots \ $ $=$
$(v_{j_1}, 2 \un{j_1} 2) \,\ldots \ $, \ where
$p_2 = (u_{j_2}, \un{j_2} 2)$.

\smallskip

\noindent Now $\Pi Q = q_1q_2q_3 \ldots \ $, for some
$q_3 = (v_{j_2}, \un{j_2} \nu) \in \Gamma$, with $\nu \in \{2,3\}$, and
$(u_{j_2},v_{j_2})$  $\in$ {\sc Input}(PCP).
Then the relation $\Pi P = \Pi Q$ takes the form
 \ $(u_{j_1} u_{j_2},\, 2 \un{j_1} 2 \un{j_2}) \,\ldots\,$ $=$
   $(v_{j_1} v_{j_2},\, 2 \un{j_1} 2 \un{j_2} \nu) \,\ldots \ \ $.

If we had $\nu = 3$ then there would be no possible choice for the
generator $p_3$ that could match the letter 3 in the second coordinate of
$q_1 q_2 q_3 \ldots \ $.
Hence, we must have $\nu = 2$. Now the relation $\Pi P = \Pi Q$ takes the
form

\smallskip

 \ \ \ $(u_{j_1} u_{j_2},\, 2 \un{j_1} 2 \un{j_2}) \,\ldots \ $ $=$
   $(v_{j_1} v_{j_2},\, 2 \un{j_1} 2 \un{j_2} 2) \,\ldots \ $,

\smallskip

\noindent where $(u_{j_1}, v_{j_1}), (u_{j_2}, v_{j_2}) \in$
{\sc Input}(PCP), and $u_{j_1} u_{j_2}$ is prefix-comparable with
$v_{j_1} v_{j_2}$.

Since $2 \un{j_1} 2 \un{j_2} \ne 2 \un{j_1} 2 \un{j_2} 2$,
$\,\Pi P$ must be of the form $p_1 p_2 p_3 \ldots \ $ for some
$p_3 \in \Gamma$ of the form $p_3 = (\e,2)$,
or $p_3 = (u_{j_3}, 2 \un{j_3})$, or
$p_3 = (u_{j_3}, 2 \un{j_3} 3)$ \ (just as at the beginning of the proof
of Claim 2).
For $p_3 = (\e,2)$ we go to case 1, and we obtain a 2-2 block.
For $p_3 = (u_{j_3}, 2 \un{j_3} 3)$ we go to case 2, and we obtain a 2-3
block. For $p_3 = (u_{j_3}, 2 \un{j_3})$, we go back to the beginning of
case 3, and the relation $\Pi P = \Pi Q$ takes the form

\smallskip

 \ \ \ $(u_{j_1} u_{j_2} u_{j_3},\,2\un{j_1}2\un{j_2}2\un{j_3})\,\ldots \ \ $
 $\,=\,$
$(v_{j_1} v_{j_2} v_{j_3},\,2 \un{j_1} 2\un{j_2} 2 \un{j_3} 2)\,\ldots \ \ $,

\smallskip

\noindent where $(u_{j_1},v_{j_1}), (u_{j_2},v_{j_2}), (u_{j_3},v_{j_3})$
$\in$ {\sc Input}(PCP), and $u_{j_1} u_{j_2} u_{j_3} $ is prefix-comparable
with $v_{j_1} v_{j_2} v_{j_3}$.
In this way case 3 could repeat itself a number of times, but since $P$ and
$Q$ have finite length, case 3 must eventually lead to case 1 or case 2;
i.e., $P$ and $Q$ both start with a 2-2 block, or both start with a 2-3
block.

\smallskip

[This proves Claim 2.]

\bigskip

\noindent {\bf Claim 3:} The relation
$\Pi P = \Pi Q$ can be factored into blocks as $ \ \Pi P = B_1\,\ldots\,B_b$
$ \ = \ $ $C_1\,\ldots\,C_b = \Pi Q$, with the following properties:

\smallskip

\noindent (3.1) \ Each block $B_i$, $C_i$ is a product of generators in
$\Gamma\,$  (for $i = 1,\ldots,b$).

\smallskip

\noindent (3.2) \ There is a matching between the generators in each pair of
blocks $(B_i, C_i)\,$ (for $i = 1,\ldots,b$); i.e., $B_i$ and $C_i$ are of
the form:

\smallskip

 \ \ \ $B_i = (u_{h_1} u_{h_2} \,\ldots\, u_{h_s},$
$\, 2\un{h_1} 2\un{h_2} \,\ldots\, 2\un{h_s} \nu)$, \ and

\smallskip

 \ \ \ $C_i = (v_{h_1} v_{h_2} \,\ldots\, v_{h_s},$
$\, 2\un{h_1} 2\un{h_2} \,\ldots\, 2\un{h_s} \nu)$,

\smallskip

\noindent where $\,(u_{h_1}, v_{h_1}), (u_{h_2}, v_{h_2}),$  $\,\ldots\, ,$
$(u_{h_s}, v_{h_s})$  $\in$ {\sc Input}(PCP), and
$\nu \in \{2,3\}$.

\smallskip

Thus, $B_i$ and $C_i$ are either both 2-2 blocks, or both 2-3 blocks.

\medskip

\noindent Notation for (3.3):
For blocks $B_i$, $C_i$, the 1st coordinates are
$\,(B_i)_1 = u_{h_1} u_{h_2} \,\ldots\, u_{h_s}$,  respectively
$\,(C_i)_1 = v_{h_1} v_{h_2} \,\ldots\, v_{h_s}$.
Similarly, the 2nd coordinates are
$\,(B_i)_2 = (C_i)_2 = 2\un{h_1} 2\un{h_2} \,\ldots\, 2\un{h_s} \nu$.

\smallskip

\noindent (3.3) \ For every $j = 1,\ldots,b:$

\smallskip

 \ \ \ $(B_1)_1\, (B_2)_1 \ \ldots \ (B_j)_1 \ $ is prefix-comparable with
$ \ (C_1)_1\, (C_2)_1 \ \ldots \ (C_j)_1$.

\smallskip

Moreover,

\smallskip

 \ \ \ $(B_1)_1\, (B_2)_1 \ \ldots \ (B_b)_1$  $=$
$(C_1)_1\, (C_2)_1 \ \ldots \ (C_b)_1$;

\smallskip

and this common string is a solution of the PCP.

\medskip

\noindent {\sf Remark:} The numbers of blocks in $P$ and $Q$ are the same
(that number is called $b$ above);
but the lengths (over $\Gamma$) of $B_i$ and $C_i$ (and hence, of $P$ and
$Q$) can be different; it depends on the number of blocks that start or end
with $(\e,2)$.

For example, the second relation in the $[\Leftarrow]$-part of this
proof has two blocks on either side.

\medskip

\noindent Proof of Claim 3: By Claim 2, $P$ and $Q$ both start with a
2-2 block, or they both start with a 2-3 block.

It could happen that $P$ consists of one such block, which implies that
$Q$ also has only one block.
Indeed, the left-most block of $Q$ has the same 2nd coordinate as the
left-most block of $P$, so the equality $\Pi P = \Pi Q$ would be violated
in the 2nd coordinate if $Q$ had additional blocks.
Similarly, if $Q$ consists of one block, then $P$ has only one block.
Since $\Pi P = \Pi Q$ holds, the 1st coordinate now yields
$\,u_{j_1} \,\ldots\, u_{j_r}$  $=$
$v_{j_1} \,\ldots\, v_{j_r}\,$ as a solution of the PCP.

If $P$ and $Q$ do not consist of just one block we look at the generators
to the right of the first block in $P$ and in $Q$. The relation
$\Pi P = \Pi Q$ then has the form
$\Pi P = B_1\,p_{s+1} \,\ldots \ = C_1\,q_{t+1} \,\ldots \ = \Pi Q$, where
$s$ and $t$ are the length of $B_1$, respectively $C_1$, over $\Gamma$.
Just as in Claim~1, the next generator after the left-most block is $(\e,2)$
on one side of the relation, and an element of
$\{(u_i, 2 \un{i}),\, (u_i, 2 \un{i} 3) : 0 \le i < k\}$ on the other side
of the relation.
As in Claim 2, this will produce either two 2-2 blocks or two 2-3 blocks just
to the right of $B_1$ and $C_1$; the process is the same as the construction
of $B_1$ and $C_1$.

Eventually, $P$ and $Q$ are factored into a finite number of blocks. The
equality $\Pi P = \Pi Q$ implies equality in the 1st and 2nd coordinates.
The 2nd coordinate gives the sequence of pairs chosen in {\sc Input}(PCP).
The 1st coordinate gives a solution of the PCP.

[This proves Claim 3.]

\medskip

\noindent By Claim 3, any non-trivial relation in
$\langle \Gamma \rangle$ determines a solution of the PCP in the 1st coordinate.  \  This concludes the proof of Lemma \ref{LEMreduction}.
\hfill $\Box$

\begin{lem} \label{LEMreduction2} {\bf (reduction without letter 3).}
 \ The symPCP with input $\big((u_i,v_i) : 0 \le i < k\big)$ over $\{0,1\}$
has a solution if and only if the subsemigroup of  $\,\{0,1\}^* \x \{0,1,2\}^*\,$
generated by

\medskip

 \ \ \  $\{(\e,2)\}$  $ \ \cup \ $
$\{(u_i, 2 \un{i}),\, (v_i, \un{i} 2)$  $\,:\, 0 \le i < k\}$

\medskip

\noindent is not free over the given generators. This holds if and only if the subsemigroup of $\,{\mathbb N}^{3 \times 3}\,$ generated by

\medskip

 \ \ \ $\{L\} \,\cup\, \{U_i,\, V_i \,:\, 0 \le i < k\}$

\medskip

\noindent is not free over the given generators.
\end{lem}
{\sc Proof.} This follows from the proof of Lemma \ref{LEMreduction} and
Remark (\ref{LEMreduction}R).
\hfill $\Box$

\bigskip

\noindent As a consequence of Lemma \ref{LEMreduction} we obtain the
following undecidability results about the freeness problem (which was
defined at the beginning of Section 3).

\begin{pro} \label{PROPundecid}
 \ The freeness problem of finitely generated subsemigroups is undecidable for the semigroups and groups below:

\smallskip

the matrix monoid $\,{\mathbb N}^{3 \times 3}$;

\smallskip

the direct product of free monoids $\,\{0,1\}^* \x \{0,1\}^*$;

\smallskip

the direct product of free groups $\,{\rm FG}_2 \x {\rm FG}_2$.

\smallskip

\noindent The number of generators of the subsemigroup can be kept bounded
without changing undecidability.
\end{pro}
{\sc Proof.} Lemma \ref{LEMreduction} yields this for
$\{0,1\}^* \x \{0,1,2,3\}^*$, and by Lemma \ref{LEMembtodirect} it then
follows for ${\mathbb N}^{3 \times 3}$.
It follows for $\{0,1\}^* \x \{0,1\}^*$, since $\{0,1,2,3\}^*$ is embeddable
into $\{0,1\}^*\,$ (e.g., by coding $\{0,1,2,3\}$ to
$\{00, 01, 10, 11\}$).  It then follows for ${\rm FG}_2 \x {\rm FG}_2$
since $\{0,1\}^*$ is a submonoid of the two-generator free group
${\rm FG}_2$.
Boundedness of the number of generators follows from Proposition~\ref{PROPsymPCPundec}.
\hfill $\Box$

\medskip

\noindent Boundedness of the number of generators was already observed in
\cite{CaHaKa}.

\bigskip

\noindent It remains an open problem whether the freeness problem for
finitely generated sub{\sl groups} of ${\rm FG}_2 \x {\rm FG}_2$ is
undecidable.

\bigskip

\bigskip




\bigskip

\bigskip

{\small
\noindent Jean-Camille Birget \\
Rutgers University - Camden \ (Emeritus) \\
Camden, NJ 08102, USA
 \ -- \ {\tt birget@camden.rutgers.edu}

\medskip

\noindent Alexey L.\ Talambutsa \\
Steklov Mathematical Institute of the Russian Academy of Sciences \\
Gubkina Str.\ 8, 119991, Moscow, Russia
  \ -- \ {\tt altal@mi-ras.ru}
}

\end{document}